\documentstyle{amsppt}
\pageheight{7.8in}
\magnification=\magstep1
\hyphenation{co-deter-min-ant co-deter-min-ants pa-ra-met-rised
pre-print
fel-low-ship}
\def\leaderfill{\leaders\hbox to 1em{\hss.\hss}\hfill}
\def\A{{\Cal A}}

\def\ldescent#1{{\Cal L (#1)}}
\def\rdescent#1{{\Cal R (#1)}}
\def\lcell#1{{{\bold L}(#1)}}
\def\rcell#1{{{\bold R}(#1)}}
\def\tcell#1{{{\bold T}(#1)}}
\def\invol#1{{\iota(#1)}}

\def\idest{i.e.\ }

\def\supp{\text{\rm supp}}
\def\boxit#1{\vbox{\hrule\hbox{\vrule \kern3pt
\vbox{\kern3pt\hbox{#1}\kern3pt}\kern3pt\vrule}\hrule}}
\def\rabbit{\vbox{\hbox{\kern0pt
\vbox{\kern0pt{\hbox{---}}\kern3.5pt}}}}

\def\tla#1{TL(\widehat{A}_{#1})}

\font\jeremy=msbm10
\def\subsetneq{\mathrel{\hbox{\jeremy \char"24}}}

\def\xfana{{\bf 1}}
\def\xfanb{{\bf 2}}
\def\xrmgj{{\bf 3}}
\def\xjona{{\bf 4}}
\def\xkl{{\bf 5}}
\def\xluse{{\bf 6}}
\def\xms{{\bf 7}}
\def\xmsa{{\bf 8}}
\def\xshi{{\bf 9}}
\def\xwes{{\bf 10}}

\topmatter
\title On the affine Temperley--Lieb algebras
\endtitle

\author C.K. Fan and R.M. Green\endauthor
\affil 
Institute for Advanced Study\\
Princeton, NJ 08540\\ USA\\
{\it E-mail:} ckfan\@math.harvard.edu\\  
\\
Mathematical Institute\\ Oxford University\\ 24--29 St. Giles'\\
Oxford OX1 3LB\\ England\\ 
{\it  E-mail:} greenr\@maths.ox.ac.uk
\endaffil

\abstract
We describe the cell structure of the affine Temperley--Lieb algebra with
respect to a monomial basis.  We construct
a diagram calculus for this algebra.
\endabstract

\thanks
The first author was supported in part by a National Science Foundation
postdoctoral fellowship.
The second author was supported in part by an E.P.S.R.C. postdoctoral
research assistantship.
\endthanks
\endtopmatter

\centerline{\bf To appear in the Journal of the London Mathematical Society}

\head 1. Introduction \endhead

Let $\Gamma$ be a Coxeter graph of type $\widehat{A}_l$, with $l>1$.
Let $W$ be the associated Coxeter group with generating set $S$.  Recall
that $S$ consists of $n=l+1$ generators, one for
each node in $\Gamma$.  These generators
satisfy the relations $s^2=1$, $st=ts$ if $s$ and $t$ are not connected in
$\Gamma$,
and $sts=tst$ if $s$ and $t$ are connected in $\Gamma$.

Every $w \in W$ may be written as a product $s_1s_2s_3\cdots s_m$ of generators
in $S$.  If $m$ is minimal, we call this product ``reduced'' and
define $\ell(w)=m$.  More generally, if $w=w_1w_2w_3 \cdots w_m$
satisfies $\ell(w)=\sum_i \ell(w_i)$, then we call this product
``reduced'' as well.

Let $H$ be the Iwahori-Hecke algebra associated to $W$.  This is
an algebra over $\A = {\Bbb Z}[v, v^{-1}]$
with generators $T_s$ for each $s \in S$ and
relations $T_s^2=(v^2-1)T_s + v^2$, $T_sT_t=T_tT_s$ if $st=ts$, and
$T_sT_tT_s=T_tT_sT_t$ if $sts=tst$, where $s$, $t \in S$.
This algebra has a basis $T_w$, $w \in W$,
where we have $T_w = T_{s_1} \cdots T_{s_m}$ whenever $s_1\cdots s_m$
is a reduced expression for $w$.

Let $I$ be the two-sided ideal generated by the elements
$$T_{sts}+T_{st}+T_{ts}+T_s+T_t+1$$
where we have one such expression for each pair of non-commuting generators
$s$, $t \in S$.

Define $\tla{l} = H / I$. 

The algebra $\tla{l}$ is infinite dimensional and is 
the affine version of the Temperley--Lieb algebra.  It has already
appeared, in the case where
$l$ is odd, in the physics literature ([\xms], [\xmsa]).
The purpose of this paper is to understand
the structure of this algebra.

We shall do this in two ways.  In the first part of the paper, 
we shall construct cells with respect
to a basis of monomials and analyze the cell structure.  This method of
attack has its roots in [\xkl] and [\xluse], but it should be
noted that the approach we adopt here is much more elementary and does
not involve Kazhdan--Lusztig theory.

In the second part of this paper, we shall
give an interpretation of $\tla{l}$ as an algebra of certain diagrams
about a cylinder.  These diagrams turn out to give the same basis as
the monomial basis mentioned above.
This is a generalization of the diagram
calculus for the ordinary Temperley--Lieb algebra, described
for example in [\xwes].  It is also a generalization of the diagram
calculus for Jones's annular algebra (see [\xjona]), which is a certain
finite dimensional algebra related to the Hecke algebra of type
$\widehat A$.  (The relationship between the affine Temperley--Lieb
algebra and the Jones algebra is explored in depth in [\xrmgj].)
The idea of representing the
affine Temperley--Lieb algebra by diagrams inscribed on a cylinder
has
been used in [\xms], but that paper does not prove that the
representation is faithful.  This is the main result of
the second part of the paper, Theorem 4.5.2.  

\vskip 20pt

\head 2. Preliminaries \endhead

\subhead 2.1 Some notation \endsubhead

Let $w \in W$.  Define the support of $w$ to be the set
of generators which appear in some (or any) reduced word for $w$.
We shall also denote this by $\supp(w)$.
Define $\ldescent{w}$ to be the set of $s \in S$
such that $\ell(sw) < \ell(w)$.  Similarly, define $\rdescent{w}$
to be the set of $s \in S$
such that $\ell(ws) < \ell(w)$.

Denote by $W_c$ those elements of $W$ whose
reduced expressions avoid substrings of the form $sts$ where $s$ and $t$
are non-commuting generators in $S$.

Elements of $W_c$ enjoy the following property, proven in [\xfana].

\proclaim{Property R}
Let $w \in W_c$, $t\in S$.  Assume $wt \not \in W_c$.  Then we may write
$w=w_1tsw_2$ reduced, $t$ commutes with all
$u \in \supp(w_2)$, and $(st)^3=1$.  Furthermore, if $w=w'_1ts'w'_2$
is another such expression, then $s=s'$.
\endproclaim

We shall refer to the dual statement of property R, which
involves generators being written in opposite order, also as property R.

For each simple generator $s \in S$, let $E_s$ be the projection
of $v^{-1}(T_s+1) \in H$ to $\tla{l}$.  For any $w \in W_c$, it makes
sense to define $E_w = E_{s_1}\cdots E_{s_n}$, where $s_1\cdots s_n$
is any reduced expression for $w$.  It was shown in [\xfana] that
the $E_w$, $w \in W_c$, constitute a basis for $\tla{l}$, and that
$\tla{l}$ is generated by the $E_s$, $s \in S$, with the
following relations:
$E_s^2=[2]E_s$, $E_sE_t=E_tE_s$ if $st=ts$, and $E_sE_tE_s=E_s$ if $sts=tst$,
$s \ne t$,
where $[2]=v+v^{-1}$.

\definition{Definition 2.1.1}
Let $w \in W_c$.  Let $s \in \ldescent{w}$.  Suppose that there
exists $t \in S$ such that $sts=tst$ and $E_tE_w=E_{sw}$.  We shall
then say that $s$ can be left cancelled, or that $s$ is left cancellable, and
that $s$ is left cancellable by $t$.
\enddefinition
We make an analogous definition for $s \in \rdescent{w}$.

\example{Example 2.1.2}
Consider $\tla{4}$.  Label the simple generators $s_1$, \dots, $s_5$,
where $s_i$ and $s_j$ commute if and only if $1<|i-j|<4$.  In the
element $w=s_1s_3s_2s_4$, the generator $s_3$ is in $\ldescent{w}$
and is left cancellable, for $E_{s_4}E_w = E_{s_1 s_2 s_4}$.
Similarly, $s_2$ is right cancellable by $s_1$.
\endexample

Let $P$ denote
the set of subsets of $\Gamma$ which consist of non-adjacent
nodes.  We allow $P$ to include the empty set, $\emptyset$.
For any $U \in P$, let $\invol{U}$ be the product of
the simple generators corresponding to the nodes in $U$ (with
$\invol{\emptyset}=1$).  Note that
the order of the product is immaterial since the nodes in $U$ are
non-adjacent.

\vskip 20pt

\subhead 2.2 Cells \endsubhead

We now imitate the definitions of left, right, and two-sided cells
given by Kazhdan and Lusztig in [\xkl].  The following definitions
have already been used in [\xfanb].  Because of the simplicity
of the structure constants for the basis $\{ E_w \}$, the definitions
of the various cells is simpler than that in [\xkl].  We remark that
these cells do not agree in general with the Kazhdan-Lusztig cells for
the affine Weyl group $W$ which are described in [\xshi].

For any $w$, $w' \in W_c$, we say that $w' \le_L w$ if there exists
$E_x$ such that $E_xE_w = c E_{w'}$ where
$c \ne 0$.

For any $w$, $w' \in W_c$, we say that $w' \le_R w$ if there exists
$E_x$ such that $E_wE_x = c E_{w'}$ where
$c \ne 0$.

For any $w$, $w' \in W_c$, we say that $w' \le_{LR} w$ if there exist
$E_x$, $E_y$ such that $E_xE_wE_y = c E_{w'}$ where
$c \ne 0$.

We say $w \sim_L w'$ if both $w' \le_L w$ and $w \le_L w'$.  Similarly,
we define $w \sim_R w'$ and $w \sim_{LR} w'$.

The equivalence classes of $\sim_L$, $\sim_R$, and $\sim_{LR}$ are
respectively, the left, right, and two-sided cells of $W_c$.
Define $\lcell{w} = \{ w' \in W_c \mid w' \sim_L w \}$,
$\rcell{w} = \{ w' \in W_c \mid w' \sim_R w \}$, and
$\tcell{w} = \{ w' \in W_c \mid w' \sim_{LR} w \}$.

\vskip 20pt

\subhead 2.3 The $a$-function \endsubhead

The following function plays the analogous role in our situation as
the Lusztig $a$-function (see [\xluse]) does for Kazhdan-Lusztig cells.

For any $w \in W_c$, we define a nonnegative integer $a(w)$ as follows.

\definition{Definition 2.3.1}
Let $w \in W_c$.
Define $$a(w)=\max _{U \in P} \{ \#U \mid
\text{\rm $\invol{U}$ is a substring of some reduced expression for
$w$}\}.$$
\enddefinition

Observe that if $\ell(sw)>\ell(w)$, then $a(sw) \ge a(w)$ since
any substring of a reduced expression for $w$ would be a substring for
a reduced expression for $sw$.  Also note that $a(\invol{U})=\# U$.

\vskip 20pt

\subhead 2.4 Left and right decompositions \endsubhead

Let $w \in W_c$.

Note that we can uniquely
write $w=\invol{G_1}\cdots \invol{G_m}$ reduced, where
$G_k=\ldescent{\invol{G_{k-1}}\cdots \invol{G_1}w} \in P$
and $G_m \ne \emptyset$.  Associated to this decomposition, we can
define graphs $\gamma_k$, $1 \le k < m$, where $\gamma_k$ is the subgraph
of $\Gamma$ which includes the nodes in $G_k$ and $G_{k+1}$ along with
all edges whose endpoints are in $G_k \cup G_{k+1}$.
We shall refer to the sequence $(G_k, \gamma_k)$ as the left
decomposition of $w$.

Similarly, we can define the right decomposition.  All the
statements we make about the left decomposition have a natural
counterpart for the right decomposition.

\proclaim{Lemma 2.4.1}
Let $w \in W_c$ be such that $w=x\invol{G_1}\invol{G_2}\invol{G_3}x'$ reduced.
We have $s \in G_1 \cap G_3$ if and only
if $s$ fails to commute with $2$ generators in $G_2$.  That is, $G_2$ contains
both generators which do not commute with $s$.
\endproclaim

\demo{Proof}
Observe that $s\invol{G_2}s$ occurs as a reduced subword of some reduced
expression for $w$.  If the lemma is false, then either $w$ is not reduced
or it is not in $W_c$.
\qed\enddemo

\vskip 20pt

\head 3. The structure of cells \endhead

\subhead 3.1 The set $Q$ \endsubhead

Let $Q \subset W_c$ denote those elements $w$ of $W_c$ such that no element
of $\ldescent{w}$ is left cancellable and no element of $\rdescent{w}$ is
right cancellable.  This leads to the following lemma.

\proclaim{Lemma 3.1.1}
Every two-sided cell contains an element of $Q$.
\endproclaim

\demo{Proof}
Let $T$ be any two-sided cell, and pick any $w \in T$.  If $w \not \in Q$,
then we can cancel a simple generator $s \in S$ either from the left or
from the right of $w$.
We continue left or right cancelling until we arrive at an element of $Q$.
The statement follows.
\qed\enddemo

When $n$ is even, there are two elements of $P$ (see 2.1) of
maximal cardinality.  Let $M_1$ and $M_2$ denote these maximal elements of
$P$.  In this case, denote by $M$ the set of elements of the form
$(\iota(M_a)\iota(M_b))^k\iota(M_a)$ and
$(\iota(M_a)\iota(M_b))^{k+1}$, for $k \ge 0$ and
positive integers $a$, $b$ such that $ab=2$.
Note that $M$ is nonempty only when $n$ is even.

\proclaim{Proposition 3.1.2}
The set $Q$ consists of the elements $\invol{T}$, $T \in P$ and, additionally,
in the case where $n$ is even, the elements of $M$.
\endproclaim

\demo{Proof}
First note that all the elements listed in the statement of the proposition
belong to $Q$.

Let $w \in Q$.  Let $(G_k,\gamma_k)$, $1 \le k \le m$,
be the left decomposition of $w$.

If $m=1$, there is nothing to prove, so assume $m>1$.

Because no element of $\ldescent{w}$ can be left cancelled, all the
nodes in $G_2$ must be connected to two nodes of $G_1$.
There are two cases.  Either $\gamma_1$ is a union of graphs of type $A$,
or $\gamma_1$ is a graph of type $\widehat{A}$.

Suppose $\gamma_1$ is a union of graphs of type $A$.  We must
have that each connected component of $\gamma_1$ has one more node
in $G_1$ than in $G_2$.  By
lemma 2.4.1, we may conclude that $\gamma_2$ is also a union of graphs
of type $A$, each component of which has one more node in $G_2$ than
in $G_3$.  Continued application of lemma 2.4.1 shows that $\gamma_{m-1}$
is a union of graphs
of type $A$, each component of which has one more node in $G_{m-1}$ than
in $G_m$.  However, this implies that an element of $G_m \subset \rdescent{w}$
is right cancellable contrary to the assumption that $w \in Q$.

Suppose $\gamma_1$ is of type $\widehat{A}$.  By the argument
used to eliminate the case where $\gamma_1$ is a union of graphs of type $A$,
we can conclude that $\gamma_k=\Gamma$ for all $1\le k < m$.  It follows
that each $G_k$ is either $M_1$ or $M_2$, and, since the left decomposition
is reduced, we have $G_k \ne G_{k+1}$.

The proposition follows.
\qed\enddemo

We now define an equivalence relation $\sim$ on $Q$.  We begin with a lemma.

\proclaim{Lemma 3.1.3}
Let $q$, $q' \in Q$.  Suppose there exists $s\in S$ such that
$E_q = E_sE_{q'}E_s$.  Then there exists $t\in S$ such that
$E_{q'}=E_t E_q E_t$.
\endproclaim

\demo{Proof}
Note that if $n$ is even and $q' \in M$, then there is no $s\in S$ such that
$E_q=E_sE_{q'}E_s$.

Therefore $q'=\invol{T}$ for some $T \in P$ such that $\#T < n/2$.
Since $E_q=E_sE_{q'}E_s$, we cannot have $s \in T$ and we cannot
have that $s$ commutes with all generators in $T$.  If $s$ fails to
commute with two generators in $T$, then $q=sq's$ reduced, which
contradicts proposition~{3.1.2}.  Therefore, $s$ fails to commute
with exactly one generator $t \in T$ and
$q=\invol{T \cup \{s\} \setminus \{t\}}$.

We compute that $E_{q'}=E_t E_q E_t$.  The lemma follows.
\qed\enddemo

We say that
$q$ and $q' \in Q$ are neighbours if and only if there exists $s\in S$ such
that $E_q=E_sE_{q'}E_s$.  The previous
lemma shows that this is symmetric.  Define
$\sim$ to be the transitive closure of the relation of being neighbours.

\proclaim{Proposition 3.1.4}
The equivalence classes of $\sim$ are as follows.

\item{\rm 1.}{$\{ \invol{T} \mid T \in P, \#T = k\}$ where $k$
is a fixed integer less than $n/2$.}

\item{\rm 2.}{The one element sets $\{ q \}$ for $q \in M$ (if $n$ is even).}
\endproclaim

We omit the proof.

\vskip 20pt

\subhead 3.2 Multiplication \endsubhead

We prove a number of facts about the algebra multiplication.

\proclaim{Lemma 3.2.1}
Let $q \in Q$.  Let $s_1, s_2, s_3, \dots, s_a$ and
$t_1, t_2, t_3, \dots, t_b$ be elements of $S$.  Define $w$ by
$E_{s_a} \cdots E_{s_1} E_q E_{t_1} \cdots E_{t_b}= cE_w$, $c \in \A$.
Then there exists $q' \in Q$ such that $q' \sim q$ and
some reduced expression for $w$ involves a substring equal to $q'$.
\endproclaim

\demo{Proof}{Let $w'=xqy$ be reduced and $s \in S$.  Define $w''$
by $cE_{w''}=E_sE_{w'}$.}

If $sw' \in W_c$ then either $w''=w'$ or $w''=sw'$ depending on whether
$\ell(sw')<\ell(w')$ or $\ell(sw') > \ell(w')$, respectively.
In either case, $w''$ has a reduced expression with a substring equal to $q$.

Assume $sw' \not \in W_c$.  Using property R, we argue as follows.

If $s$ commutes with all generators in $\supp(xq)$, then
$E_sE_{w'}=E_xE_qE_sE_y$.  Note that $E_sE_y$ may be written
as a product of $\ell(y)-1$ basis elements $E_{s_k}$ for $s_k \in S$.

If $s$ commutes with all generators in $\supp(x)$,
but not with all the generators in $\supp(q)$, then
we must have $y=sy'$ and $s$ commutes with all but one $t \in \supp(q)$.
Let $T=(\supp(q)\cup \{s\} )\setminus \{t\}$.  Observe that $\invol{T} \sim q$.
We compute that
$E_sE_{w'}=E_xE_{\invol{T}}E_{y'}$.

If $s$ does not commute with all the generators in $\supp(x)$, then
$E_sE_x$ can be written as a product of $\ell(x)-1$ basis elements
$E_{s_k}$ for $s_k \in S$.

By symmetry and induction on the length of the sequence $s_k$,
these observations imply the lemma.
\qed\enddemo

We have the following corollary of lemma~{3.2.1}.

\proclaim{Corollary 3.2.2}
Let $w \in W_c$, $s \in S$.  Define $w'$ by
$c E_{w'}=E_wE_s$, $c \in \A$.  Then $a(w') \ge a(w)$.
\endproclaim

Lemma~{3.2.1} can be strengthened in the following situation.
\proclaim{Lemma 3.2.3}
Let $q \in Q$.  Let $s_1$, $s_2$, $s_3$, \dots, $s_m$ be elements of $S$.
Define $w' \in W_c$ by
$E_q E_{s_1}E_{s_2}E_{s_3}\cdots E_{s_m}=cE_{w'}$, $c \in \A$.
Then one can write $w'=qx$ reduced.
\endproclaim

\demo{Proof}{We proceed by induction on $m$, the case where $m=0$ being clear.}

Suppose the lemma is true for all $m<N$.  We prove the result holds for
$m=N$.

By induction, there exists $y$ such that
$E_q E_{s_1}E_{s_2}E_{s_3}\cdots E_{s_{m-1}}=c'E_{y}$ and
$y=qy'$ reduced.  Define $w'$ by $E_{w'}=cE_{y}E_{s_m}$.

There are three cases to consider.  Either (1) $\ell(ys_m)<\ell(y)$,
(2) $\ell(ys_m)>\ell(y)$ and $ys_m \in W_c$, or
(3) $ys_m \not \in W_c$.

In case 1, $c=[2]$ and $w'=y$, and the lemma holds.

In case 2, $c=1$ and $w'=ys_m$, and the lemma holds.

In case 3, let $w_1s_mtw_2$ be a decomposition of $y$ as in property
R.  We may further assume that $w_1s_m=qx'$, reduced, by commuting
to the left of $s_mt$ all possible $s \in \supp(w_2)$.
We have $E_{w_1}E_{s_m}E_tE_{w_2}E_{s_m} = E_{w_1}E_{s_m}E_{w_2}$.
Although the product $w_1s_mw_2$ is not necessarily reduced, it represents
a product of generators of length less than $N+\ell(q)$.  Thus, the lemma
follows by induction.
\qed\enddemo

\proclaim{Lemma 3.2.4}
Let $w=xq$, reduced, with $q \in Q$ and $a(w)=a(q)$.
Assume that no element of $\rdescent{w}$
is right cancellable.  Let $s_1, s_2, s_3, \dots, s_j$ be elements of
$S$.  Define
$w' \in W_c$ by $c E_{w'} = E_w E_{s_1} E_{s_2} E_{s_3} \cdots E_{s_j}$,
$c \in \A$.
If $a(w')=a(q)$, then there
exists $y \in W_c$ such that $w'=xqy$ reduced.
\endproclaim

\demo{Proof}
We prove this by induction on $j$, the case where $j=0$ being
clear.

Suppose the lemma is true for $j<N$.  Let $z\in W_c$ be defined by
$c' E_z = E_w E_{s_1} E_{s_2} E_{s_3} \cdots E_{s_{N-1}}$, 
$c' \in \A$.  If $a(z) \ne a(q)$,
then $a(z)>a(q)$ and by corollary~{3.2.2}, this would imply
that $a(w')>a(q)$, and there would be nothing to prove.

So assume that $a(z)=a(q)$.  Then by induction, there
exists $y' \in W_c$ such that $z=x\invol{U}y'$ reduced.  Define $u$ by
$d E_u = E_zE_{s_N}$, $d \in \A$.
It suffices to show that $u = x\invol{U}y''$
reduced, for some $y'' \in W_c$.

To see this, we consider the standard cases.  If $zs_N \in W_c$, then
$u=z$ or $u=zs_N$ depending on whether $\ell(zs_N)<\ell(z)$ or
$\ell(zs_N)>\ell(z)$, respectively.  In either case, there is nothing to show.

So assume $zs_N \not \in W_c$.  In this case, write $z=x_1s_Ntx_3$
as in property~R.  Observe that $t \not \in \rdescent{q}$
since by assumption,
no element of $\rdescent{w}=\rdescent{q}$
is right cancellable.  On the other hand,
if $t \in \supp(x)$ and $t \not \in \supp(y')$, we must have, by
property~R and maximality of $\# U$, 
that $s_N$ commutes with all generators in $U$.  Consequently
$a(w')>a(q)$, and there is nothing to prove.  The only remaining possibility
is that $t \in \supp(y')$ so that
$E_zE_{s_N} = E_{x_1} E_{s_N} E_{x_3} = E_{x\invol{U}}E_{y'}E_{s_N}$,
where $E_{y'}E_{s_N}$ can
be written as a product of $\ell(y')-1$ algebra generators $E_s$.
We can therefore apply induction and the lemma follows.
\qed\enddemo

\vskip 20pt

\subhead 3.3 Two-sided cells \endsubhead

We now classify the two-sided cells $\tcell{w}$ which were defined in 2.2.

\proclaim{Theorem 3.3.1}
The two-sided cells of $W_c$ are parametrized by $Q/\sim$ (see 3.1).
The bijection is given by $q \mapsto \tcell{q}$.  Furthermore,
the $a$-function is constant on two-sided cells.
\endproclaim

\demo{Proof}
We begin with the following lemma.

\proclaim{Lemma 3.3.2}
Let $q$, $q' \in Q$.
We have $q \sim q'$ (see 3.1) if and only if $q \sim_{LR} q'$.
\endproclaim

\demo{Proof}{First we show that $q \sim q'$ implies
$q \sim_{LR} q'$.
It suffices to check this in the case where $q$ and $q'$ are neighbours.}

By definition, there exists
$s$, $t\in S$ so that $E_q=E_sE_{q'}E_s$
and $E_{q'}=E_tE_q E_t$.  This shows that $q \sim_{LR}
q'$.

For the other implication, suppose $q \not\sim q'$.  By lemma~{3.2.1},
every $w$ such that there exists a constant $c$ for which
$cE_w=E_xE_qE_y$ must have a reduced expression which involves
a substring equal to some $q''$ such that $q'' \sim q$.  If
$q \sim_{LR} q'$, then some substring of $q$ is equal to $q'$ and
some substring of $q'$ is equal to $q$.  That is, $q=q'$, contrary
to hypothesis.  Thus, $q \not \sim_{LR} q'$.
\qed\enddemo

Theorem~{3.3.1} follows from lemma~{3.3.2} and
lemma~{3.1.1}.
\qed\enddemo

\vskip 20pt

\subhead 3.4 Involutions in $W_c$ \endsubhead

Before we classify the right cells, we first study involutions
in $W_c$.

\proclaim{Theorem 3.4.1}
Let $w \in W_c$ satisfy $w^2=1$.  There exists a unique $T \in P$
and a unique $x \in W_c$ such that $w=x\invol{T}x^{-1}$ reduced.
\endproclaim

Because of the uniqueness of $x$ and $T$ in theorem 3.4.1,
we shall call $x\invol{T}x^{-1}$ the canonical decomposition of
an involution.

\demo{Proof}  The proof is the same as the proof of
theorem 3.1.1 in [\xfanb].
\qed\enddemo

Recall from the proof of theorem 3.1.1 in [\xfanb] that one can construct a
finite sequence of triples $(w_m, x_m, S_m)$
for $0 \le m \le N$, such that
$w_m$, $x_m \in W_c$, $w_m^2=1$, $w=x_mw_mx_m^{-1}$ reduced,
$S_m = \supp(w_m)$, $\ell(x_{m+1}) = \ell(x_m) + 1$, 
and $N$ is minimal such that $S_N \in P$.

In type affine $A$, we also note that
$\supp(w_m)$ is either equal to $\supp(w_{m+1})$
or differs by a single node.  Also, because no element of $W_c$ whose
support is a subdiagram of a graph of type $A$
involves two occurrences of an end generator, we conclude that
the number of connected components of $\supp(w_m)$ is exactly $1$ less than
the number of connected components of $\supp(w_{m+1})$ when
$\supp(w_{m+1}) \ne \supp(w_m)$ and $\supp(w_m) \ne \Gamma$.  It follows that
none of the graphs $\supp(w_m)$ can contain a subgraph of type $A_2$
(since the two elements of $W_c$ for type $A_2$ whose support is
all of $A_2$ are not involutions).
Thus, we have the following corollary.

\proclaim{Corollary 3.4.2}
When $n$ is odd, no involution $d \in W_c$ has $\supp(d)=\Gamma$.
If $\supp(d)=\Gamma$, then $T$ is maximal with respect to inclusion in $P$,
where $d=x\invol{T}x^{-1}$ is the canonical decomposition.
\endproclaim

Lemma~{2.4.1} implies the following.

\proclaim{Corollary 3.4.3}
Let $d \in W_c$ be an involution and let $d=x\invol{T}x^{-1}$ be the
canonical decomposition.  Then no element of $\rdescent{x\invol{T}}=T$
can be right cancelled.
\endproclaim

\proclaim{Proposition 3.4.4}
Let $d \in W_c$ be an involution and let $d=x\invol{T}x^{-1}$ be the
canonical decomposition.  Assume that $\supp(d) \ne \Gamma$.
Then $d \sim_R x\invol{T}$,
$d \sim_L \invol{T}x^{-1}$, and $d \sim_{LR} \invol{T}$.
\endproclaim

\demo{Proof}
We assume $x \ne 1$, otherwise the statement is clear.

Consider the element $x\invol{T}$.  Note that $\rdescent{x\invol{T}}=T$
since $x\invol{T}x^{-1}$ is reduced.
Note that no element of $\rdescent{x\invol{T}}$ is right cancellable
since $x\invol{T}x^{-1} \in W_c$.

By proposition~{3.1.2}, 
we conclude that some element of $\ldescent{x\invol{T}}$ is left cancellable.

Using induction on the length of $x$, we conclude that
$x\invol{T} \sim_L \invol{T}$.  A similar argument shows that
$\invol{T} \sim_R \invol{T}x^{-1}$.  Thus,
$x\invol{T} \sim_R x\invol{T}x^{-1}$ and
$\invol{T}x^{-1} \sim_L x\invol{T}x^{-1}$.  The proposition follows.
\qed\enddemo

\proclaim{Corollary 3.4.5}
Let $d \in W_c$ be an involution and let $d=x\invol{T}x^{-1}$ be the
canonical decomposition.  We have $a(d)=\#T$.
\endproclaim

\demo{Proof}
We have $a(d) \ge \#T$.  If $\#T$ is maximal, there is nothing to prove.
Otherwise, by corollary~{3.4.2}, $\supp(d) \ne \Gamma$.  By
proposition~{3.4.4} and the fact that the $a$-function is constant
on two-sided cells, we see that $a(d)=a(\invol{T})=\#T$.
\qed\enddemo

\vskip 20pt

\subhead 3.5 Right cells \endsubhead

We are now ready to begin classifying right
cells.

\proclaim{Theorem 3.5.1}
The map given by $d \mapsto \rcell{d}$,
defines an injection from the set of
involutions to the set of right cells.
If $n$ is odd, this map is surjective.
If $n$ is even, the right cells not in the image of this map are
the right cells in the two-sided cells parametrized by elements of $M$
which are not involutions.
\endproclaim

The remainder of this section is devoted to a proof of theorem~{3.5.1}.

First note that, by definition,
each two-sided cell is a disjoint union of right cells.
Let $q \in Q$.  We shall first consider the case where $q \not \in M$.

Let $w \in \tcell{q}$.  By right cancelling repeatedly, we can find
$x$ such that $w \sim_R x$ and no element of $\rdescent{x}$ can be right
cancelled.

By lemma~{3.2.1}, we can write $x=yq'z$ where $q' \sim q$ and
$a(x)=a(w)$.  Let $\invol{G_1}\cdots \invol{G_m}$ be the right decomposition
of $x$.  Note that every element of $G_{m-1}$ must be connected to
two elements of $G_m$, because no element of $\rcell{x}$ can be right
cancelled.  By lemma~{2.4.1} and induction, we see that
each element of $G_k$ is connected to two elements of $G_{k+1}$.
Since $\#\rdescent{x} \ne n/2$, we see that the cardinality of $G_k$
strictly increases with $k$.  In other words, $q'=\invol{G_m}$ and
$x=yq'$.

Observe that $yq'y^{-1}$ is reduced and an involution (compare lemma 4.3.1
of [\xfanb]).  By proposition~{3.4.4}, we have $yq'y^{-1} \sim_R w$.
Thus, every $w \in \tcell{q}$ is right equivalent to some involution.
By lemma~{3.2.4} and theorem~{3.4.1},
we deduce that every right cell in $\tcell{q}$
contains a unique involution.

Now assume that $q \in M$.  Then $n$ is even.

First, consider the element $M_1$.  By proposition~{3.1.4} and lemma~{3.2.3},
every element of $\rcell{M_1}$ is of the form $M_1y$ reduced.
If $M_1y \sim_R M_1$ then $\ldescent{y}$ must be a proper subset of $M_2$.
By a similar argument used in the case where $q \not \in M$, we conclude
that any $y$ such that $\ldescent{y} \ne M_2$ satisfies $M_1 \sim_R M_1y$.

These observations enable us to deduce that the right cells in the two-sided
cell $\tcell{q}$ are parametrized by $z$ such that $z\invol{\ldescent{q}} \in
W_c$ reduced and $\rdescent{z} \subsetneq S \setminus \ldescent{q}$.  We have
$$\rcell{zq}=\{ zqy \mid
\ldescent{y} \subsetneq S \setminus \rdescent{q}\}.$$

If $\rdescent{q}=\ldescent{q}$, then $zqz^{-1} \in W_c$ is reduced and
an involution.

The theorem follows. \qed

\vskip 20pt

\head 4. A Diagram Calculus for $\tla{l}$ \endhead

Label the elements of $S$ by $s_1, \dots, s_n$, where
$s_i$ and $s_j$ are connected in $\Gamma$ whenever $i$ and $j$ are
consecutive modulo $n$.

\subhead 4.1 Affine $n$-diagrams \endsubhead

\definition{Definition 4.1.1}
An affine $n$-diagram
consists of two infinite horizontal rows of nodes together with
edges which can be made to satisfy the following conditions:

\item{\rm (i)}
{Every node is the endpoint of exactly one edge.}
\item{\rm (ii)}
{Any edge lies entirely between the two rows of nodes.}
\item{\rm (iii)}
{If an edge does not link two nodes then it is an infinite horizontal
line which does not meet any node.}
\item{\rm (iv)}
{No two edges intersect each other.}
\item{\rm (v)}
{The diagram is invariant under shifting to the left
or to the right by $n$ nodes.}
\enddefinition

Because of the condition {\rm (v)},
one can also think of affine $n$-diagrams as
diagrams on the surface of a cylinder, or within an annulus, in a natural way.
Unless otherwise specified, we
shall henceforth regard the diagrams as diagrams on the surface of
a cylinder with $n$ nodes on top and $n$ nodes on the bottom.
Compare with Jones' (finite-dimensional) annular 
algebra in [\xjona].

An example of an affine $n$-diagram for $n = 4$ is given in figure 1.

\head Figure 1: An affine 4-diagram \endhead
\centerline{
\hbox to 3.638in{
\vbox to 0.888in{\vfill
	\includegraphics{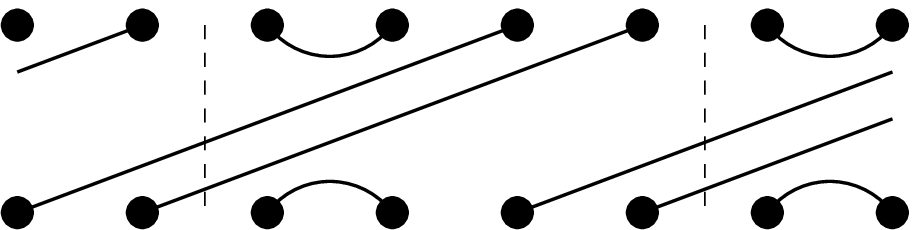}
}
\hfill}
}

\vskip 20pt

\subhead 4.2 Representation of $\tla{l}$
by affine $n$-diagrams \endsubhead

Consider the free $\A$-module with basis given by
the affine $n$-diagrams.  We construct an algebra structure on this
$\A$-module by defining a multiplication as follows.
For two diagrams $A$ and $B$,
put the cylinder for $A$ on top of the cylinder for $B$ so that corresponding
vertices match and then identify
all the points in the middle row.  This produces a certain number, say $x$,
of loops.  Removal of these loops forms another diagram $C$ satisfying
the conditions in Definition 4.1.1.  The product $AB$ is then defined
to be $[2]^x C$.
One can check that this extends to an associative multiplication with
unit.  The unit is given by the diagram for which node $\bar{j}$
in the top ring is joined to node $\bar{j}$ in the bottom
ring for all $\bar{j}$.

This definition generalizes the realization of
the finite Temperley--Lieb algebra in terms of
$r$-diagrams---see for example [\xwes, \S1].

Denote by $E'_{s_i}$ the unique diagram in which
$\bar{i}$ and $\overline{i+1}$ are joined in both rows of nodes, and
for other congruence classes $\bar{j}$, node
$\bar{j}$ in the top ring is joined to node
$\bar{j}$ in the bottom ring.


\proclaim{Proposition 4.2.1}
There is a unique $\A$-algebra homomorphism $\rho$ from $\tla{l}$
to the algebra of affine $n$-diagrams
which sends $E_{s_i}$ to $E'_{s_i}$.
\endproclaim

\demo{Proof}
One checks that the defining relations hold.
\qed\enddemo

The aim of the rest of \S4 is to prove that $\rho$ is injective and to
characterize a set of diagrams which constitute a basis for
$\rho(\tla{l})$.

\subhead 4.3 Monomials and diagrams \endsubhead

Proposition 4.2.1 allows us to make the following definition.

\definition{Definition 4.3.1}
Choose an element $w \in W_c$ of length $k$ and fix a reduced
expression $s_{i_1} \cdots s_{i_k}$ for $w$.  The stacked
representation of this reduced expression is the diagram,
denoted $E'_w$, formed by
stacking the diagrams $E'_{s_i}$
on top of each other, with $s_{i_1}$ at the top.
\enddefinition

\proclaim{Lemma 4.3.2}
For $w \in W_c$, the diagram $E'_w$ is an affine
$n$-diagram.
\endproclaim

\demo{Proof}
One has to check that the conditions of definition 4.1.1 are satisfied.
All are immediate except (iii).  For (iii), note that any edge either
links two nodes or is equivalent to a closed loop on the
cylinder.
\qed\enddemo

\remark{Remark 4.3.3}
Let $w \in W_c$.  Note that by construction,
$\rho(E_w) = E'_w$.
\endremark

\proclaim{Lemma 4.3.4}
Let $w \in W_c$.  Then we have: 
\item{\rm (i)}
{Points $\bar{i}$ and $\overline{i+1}$ in the top ring of $E'_w$ are
connected to each other by an edge of minimal length if and only if
$\ell(s_iw) = \ell(w) - 1$.}
\item{\rm (ii)}
{Points $\bar{i}$ and $\overline{i+1}$ in the bottom ring of $E'_w$ are
connected to each other by an edge of minimal length if and only if
$\ell(ws_i) = \ell(w) - 1$.}
\endproclaim

\demo{Proof}
Part (i) follows from the fact that $\ell(s_iw)<\ell(w)$ if and only if
$w=s_iw'$ reduced, and from the definition of $E'_w$.
Part (ii) follows from the fact that $\ell(ws_i)<\ell(w)$ if and only if
$w=w's_i$ reduced, and from the definition of $E'_w$.
\qed\enddemo

Note that the number of times $s_k$ appears in some reduced
expression for $w \in W_c$ is independent of the reduced expression
for $w$.  Denote this number by $\nu(k,w)$.

We shall have to consider auxiliary vertical lines on the cylinder.
A vertical line which passes between nodes $\bar{k}$ and
$\overline{k+1}$ will be called the vertical line $k+1/2$.

\proclaim{Lemma 4.3.5}
\item{\rm(i)}
{Let $k$ be an integer from $1$ to $n$ inclusive.
Then the number of intersections of $E'_w$ with the vertical line $k +1/2$
is equal to $2 \nu(k, w)$.}
\item{\rm(ii)}
{Let $w \in W_c$.  Then the path traced out by an edge in the diagram
$E'_w$ either has no segments moving to the left, or no segments
moving to the right.
}
\endproclaim

\demo{Note}
Informally, (ii) is saying that the edges in a diagram corresponding
to an element of $W_c$ do not double back on themselves.
\enddemo

\head Figure 2: A diagrammatic version of property R \endhead
\centerline{
\hbox to 0.527in{
\vbox to 1.277in{\vfill
	\includegraphics{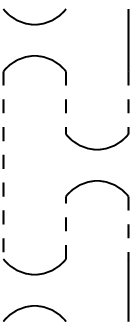}
}
\hfill}
}

\demo{Proof}
We first show (i).  
We prove this by induction on the length of $w$, the case where
$\ell(w)$ is $0$ or $1$ being clear.

Suppose the lemma holds for $w$ such that $\ell(w)<l$, where $l>1$.
Let $w \in W_c$ have length $l$.  We can write $w=w's_i$ reduced for
some generator $s_i$.  By induction, the number of intersections
of $E'_{w'}$ with the vertical line $k +1/2$
is equal to $2 \nu(k, w')$.  There are two cases: either
$i=k$ or $i \ne k$.

Suppose $i \ne k$.  Then by definition of $E'_{s_i}$ and $E'_w$, we see
that the number of intersections of $E'_w$ with the vertical line $k+1/2$ is
equal to the number of intersections of $E'_{w'}$ with the vertical line
$k+1/2$ and the lemma holds for $w$.

Suppose $i =k$.  By lemma 4.3.4, we know that nodes $\bar{k}$ and
$\overline{k+1}$ are not connected by an edge of minimal length so that
no loops are formed in the diagram multiplication of $E'_{w'}$ with
$E'_{s_k}$.  Therefore, the number of intersections of $E'_w$
with the vertical line $k+1/2$ is
equal to the number of intersections of $E'_{w'}$ with the vertical line
$k+1/2$ plus $2$, and the lemma holds for $w$.

We now prove (ii), again by induction on $\ell(w)$.  The base
cases are clear.  A short consideration shows that a counterexample
of minimal length must consist of a diagram containing a subdiagram like
the one shown in figure 2, or its left-right 
mirror image.  (The dashed lines denote lines
of indeterminate length.)  However, this shows that the
corresponding element $w$ contains a reduced subword $s w_1 t w_2 s$,
where $w_1$ and $w_2$ are words in the generators,
$s$ and $t$ are generators which
do not commute with each other, and $s$ commutes
with $w_1$ and $w_2$.  Thus $w$ has a reduced expression containing
the sequence $sts$, and does not lie in $W_c$.

The lemma follows.
\qed\enddemo

\remark{Remark 4.3.6}
The situation in figure 2 is a manifestation of Property R.
\endremark

\vskip 20pt

\subhead 4.4 Straightening diagrams \endsubhead

Consider a diagram $D$.  Perturb the edges as necessary so that edges
never run along a vertical line $k+1/2$.
Fix an edge and a node $\bar{k}$.  If the edge
has endpoints, we can compute the number of crossings of this edge with
the vertical line $k+1/2$ by counting crossings with multiplicity $+1$
or $-1$ depending on whether the crossing is anti-clockwise or
clockwise, adding these intersections multiplicities up and taking
absolute value.  If the edge has no endpoints, then we can declare the
number of crossings with the vertical line $k+1/2$ to be $+1$.
Note that these numbers are only dependent upon the graph theoretic
properties of $D$.  We define $\nu_k(D)$ to be the sum of these
(nonnegative) crossing numbers over all edges.  It follows from Lemma
4.3.5 that $\nu_k(D)=2\nu(k,w)$ when $D=E'_w$ for $w \in W_c$.

It is convenient to classify the edges in an affine $n$-diagram in the
following way.

\definition{Definition 4.4.1}
We call an edge linking two points on the same side (\idest top or
bottom) of an affine $n$-diagram a short horizontal
edge.  We call an edge with no endpoints a long horizontal edge.  We
call the other edges vertical.
\enddefinition

\definition{Definition 4.4.2}
An affine $n$-diagram $D$ is said to be admissible
if {\rm (a)} it
has no horizontal edges and is the image of the identity element in
$\tla{l}$ or {\rm (b)} it has horizontal edges and 
the number $\nu_k(D)$ of intersections of the diagram with
any vertical line $k + 1/2$ is an even number.

If $D$ is admissible, the nonnegative integer $$
{1 \over 2} \sum_{k = 1}^n \nu_k(D)
$$ is called the length of $D$, denoted by $\ell(D)$.
\enddefinition

Note that this definition implies that there are only a finite number
of long horizontal edges in any admissible diagram.  

\definition{Definition 4.4.3}
The admissible diagram $D$ is called straight if it is of form
$E'_{\invol{S}}$, $S \in P$.
\enddefinition

The next aim is to express a non-straight diagram $E'_w$ in a definite
way as a product
$E'_{s_i} E'_{w'}$ where $\ell(E'_{w'}) = \ell(E'_w) - 1$.

\definition{Definition 4.4.4}
Let $D$ be an admissible but non-straight diagram represented between
two infinite horizontal lines.

We say that a congruence class $\bar{k}$ mod $n$ is of type 1T for $D$ if
points $k$ and $k+1$ in the top row of $D$ can be connected by a short
horizontal edge of minimal length, and either 
{\rm (a)} $D$ has a short horizontal edge
connecting points $i$ and $j$ in its top row, where $i < k$ and $j >
k+1$ are integers or if {\rm (b)}
$D$ has no edges as described in (a), but has a long horizontal edge.

We define a class of type 1B as above but using the bottom edge.

We say that $\overline{k-1}$ is of type 2T for $D$ if points
$k$ and $k+1$ in the top row of $D$ can be joined to each other 
by a short horizontal
edge of minimal length and top node $k-1$ is connected to a
bottom node $j$ with $j>k+1$.

We say that a class $\overline{k-1}$ is of type 2B for $D$ if points
$k$ and $k+1$ in the bottom row of $D$ can be joined to each other
by a short horizontal
edge of minimal length and bottom node $k-1$ is connected to a
top node $j$ with $j>k+1$.
\enddefinition

\example{Example 4.4.5}
Consider the case $n = 4$.  Then 2 is a congruence class of type 1T
and 1B for the diagram $E'_w$ where $w = s_2 s_1 s_3 s_2$.  Also 1 is
a congruence class of type 2T for the diagram $E'_{w'}$ where $w' =
s_2 s_1$.
\endexample

\proclaim{Lemma 4.4.6}
Any admissible but non-straight diagram $D$ has a congruence class
$\bar{k}$ mod $n$ of type 1T, 1B, 2T or 2B.
\endproclaim

\demo{Proof}
Let us assume $D$ has no congruence classes of type 1T or 1B.

Suppose $D$ has a vertical edge.  Then not all the edges can be
vertical unless $D$ is the image of the identity element (which would
mean $D$ was straight, contrary to assumption), because
$D$ is admissible.  It now follows that $D$ will have a congruence
class of type 2T if it has
a vertical edge with negative gradient, and a congruence class of type 2B
if it has a vertical edge with positive gradient.  One of these situations
must arise under the hypotheses on $D$.  (It should be noted that no
admissible diagram can have a vertical edge of minimal length connecting a
top node $a$ to a bottom node $a \pm 1$.)

Suppose that $D$ has no vertical edges, and no congruence classes of
type 1T or 1B.  It therefore consists entirely of short horizontal
edges of minimal length.  The only such possibilities are diagrams
which are either straight, or inadmissible.
\qed\enddemo

\definition{Definition 4.4.7}
We say the congruence class $i$ mod $n$ is a distinguished congruence
for $D$ if one of the following holds:

\item{\rm (i)}
{$D$ has a congruence class $i$ of type 1T.}
\item{\rm (ii)}
{$D$ has no congruence class of type 1T, but has a congruence class
$i$ of type 1B.}
\item{\rm (iii)}
{$D$ has no congruence class of type 1T or 1B but has a congruence
class $i$ of type 2T.}
\item{\rm (iv)}
{$D$ has a no congruence class of type 1T, 1B or 2T, but has a congruence
class $i$ of type 2B.}

Note that exactly one of these conditions will hold.  We classify the
congruence to be of type (i), (ii), (iii) or (iv)
respectively depending on which condition holds.
\enddefinition

The following lemma is the main result of this section.

\proclaim{Lemma 4.4.8 (Straightening Rule)}
Let $D$ be an admissible but not straight diagram and let $i$ be a
distinguished congruence for $D$.

If $i$ is of type (i) for $D$ then there exists a unique
admissible (possibly straight) diagram $D'$ of length $\ell(D) - 1$
such that $E'_{s_i} D' = D$.

If $i$ is of type (ii) for $D$ then there exists a unique
admissible diagram $D'$ of length $\ell(D) - 1$
such that $D' E'_{s_i} = D$.

If $i$ is of type (iii) for $D$ then we define $D' =
E'_{s_i} D$.  Then $D'$ is admissible,
$\ell(D') = \ell(D) - 1$ and $D = E'_{s_{i+1}} D'$.

If $i$ is of type (iv) for $D$ then we define $D' = 
D E'_{s_i}$.  Then $D'$ is admissible,
$\ell(D') = \ell(D) - 1$ and $D = D' E'_{s_{i+1}}$.

In each case, $D'$ has the same number of short horizontal edges as $D$.
\endproclaim

\demo{Proof}
The proof of (i) is by an entirely routine exhaustive
check, so we omit it.  The proof of (ii) is essentially the same as
the proof of (i).  The proof of (iv) is analogous to the proof
of (iii), so we concentrate on proving (iii).

In the case of (iii) one verifies using the product rule that $E'_{s_i} D$ is
equal to another diagram $D'$ satisfying $\nu_{i+1}(D') =
\nu_{i+1}(D)-1$
and $\nu_k(D') = \nu_k(D)$ for $k \ne i+1$.  The claim about the
length
of $D'$ follows.  The claim that $D = E'_{s_{i+1}} D'$ is
a consequence of Lemma 4.3.4, the product rule, and
the relation $E_{s_{i+1}} E_{s_i} E_{s_{i+1}} = E_{s_{i+1}}$.
The claim about
the number of short horizontal edges also follows by direct inspection.
\qed\enddemo

\vskip 20pt

\subhead 4.5 Consequences of the Straightening Rule \endsubhead

\proclaim{Corollary 4.5.1}
Any admissible diagram is of the form $E'_w$ for some $w \in W_c$.
Conversely any $E'_w$, $w \in W_c$, is admissible.
\endproclaim

\demo{Proof}
Suppose $D$ is a counterexample to the statement of minimal length.
Then $D$ is clearly not straight.  Lemma 4.4.8 now shows $D$ to be the
product of a shorter admissible diagram and a generator $E'_{s_i}$.  By
hypothesis, $D'$ is of form $E'_{w'}$ for $\ell(w') = \ell(w) - 1$.
Since no $[2]$ appears, the proof of Lemma 4.3.4 shows that $D$ is
of form $E'_w$ where $w = s_i w'$ or $w = w' s_i$ according as $D = E'_{s_i} D'$
or $D = D' E'_{s_i}$.

By Lemma 4.3.4 and Lemma 4.3.5, one can see easily that any $E_w$ maps
under $\rho$ to an admissible diagram.
\qed\enddemo

This means the image of $\tla{l}$ in the diagram algebra is precisely
the span of the admissible diagrams.

\proclaim{Theorem 4.5.2}
The diagram representation is faithful.
\endproclaim

\demo{Proof}
Recall that a basis for $\tla{l}$ is given by
$\{ E_w \}_{w \in W_c}$.  Let $a$, $b \in W_c$.
By Remark 4.3.3, if $\rho$ were not
faithful, two monomials $E_a$ and $E_b$ would map to the same diagram
$D$.  By Lemma 4.3.4, $\rho$ preserves length, so we conclude
that $a$ and $b$ have the same length.

Now suppose we have a counterexample of minimal length.  Suppose $D$
is not straight.  Then by Lemma 4.4.8 either (i)
$D = E'_{s_i} D'$ where $D' = E'_{a'} = E'_{b'}$, $a = s_i a'$, $b = s_i b'$ 
for some $i$ or (ii)
$D = D' E'_{s_i}$ where $D' = E'_{a'} = E'_{b'}$, $a = a' s_i$, $b = b'
s_i$ for some $i$.

In either case, $E'_{a'} = E'_{b'}$ is a shorter counterexample.  So
$D$ must be straight.

By Lemma 4.3.4, if $D = E'_w$ then a reduced expression for 
$w$ must contain all the commuting
reflections associated to it by the definition of it being straight.
A quick analysis of $\nu(k, w)$ for the various $k$ using Lemma 4.3.5
now shows that $w$ is unique.  The proof follows.
\qed\enddemo

We conclude by giving the relation between cell structure and
diagrams.

\proclaim{Proposition 4.5.3}
Let $w \in W_c$.  Then $a(w)$ is equal to
the number of short horizontal edges on the top (or bottom)
face.

The two-sided cell in which $w$ lies is determined by $a(w)$ and
the number of long horizontal edges.

The left cell in which $w$ lies is determined by its two-sided cell
and the pattern of short horizontal edges on the bottom row.

The right cell in which $w$ lies is determined by its two-sided cell
and the pattern of short horizontal edges on the top row.
\endproclaim

\demo{Proof}
This follows by noting that left equivalence preserves the pattern of
short horizontal edges on the bottom row, right equivalence preserves
the pattern of short horizontal edges on the top row, and therefore,
in both cases,
the number of short horizontal edges is invariant.  Since, by
the results of section 3, for each $w \in W_c$ there exists $q \in Q$
such that $w \sim_{LR} q$, and because for $q \in Q$, we can see that
$a(q)$ is the number of short horizontal edges in $E'_q$, the
proposition follows.
\qed\enddemo

Now using standard combinatorial arguments we have the following.

\proclaim{Corollary 4.5.4}
The number of left cells (or right cells) in a two-sided cell
with $a$-value $k$, $k < n/2$, is equal to $n \choose k$.
If $k=n/2$, then the number of left cells (or right cells) in
a two-sided cell with $a$-value $k$ is ${1 \over 2}{n \choose k}$.
\endproclaim

\vskip 1cm

\head References \endhead

\item{[\xfana]}
{C.K. Fan, {\it A Hecke Algebra Quotient and Some Combinatorial Applications},
J. of Alg. Comb. {\bf 5} (1996) no. 3, 175--189.}
\item{[\xfanb]}
{C.K. Fan, {\it Structure of a Hecke Algebra Quotient}, Jour.
Amer. Math. Soc. {\bf 10} (1997) no. 1, 139--167.}
\item{[\xrmgj]}
{R.M. Green, {\it On representations of affine Temperley--Lieb
algebras}, Proceedings of the 8th International Conference on
Representations of Algebras, to appear.}
\item{[\xjona]}
{V.F.R. Jones, {\it A quotient of the affine Hecke algebra in the
Brauer algebra}, L'Enseignement Math. {\bf 40} (1994) 313--344.}
\item{[\xkl]}
{D. Kazhdan and G. Lusztig, {\it Representations of Coxeter groups
and Hecke algebras}, Invent. Math. {\bf 53} (1979) 165-184.}
\item{[\xluse]}
{G. Lusztig, {\it Cells in affine Weyl groups}, Algebraic groups and related
topics, Adv. Studies Pure Math 6, North-Holland and Kinokuniya,
Tokyo and Amsterdam, 1985, 255--287; {\it II}, J. Alg. {\bf 109} (1987),
536--548; {\it III}, J. Fac. Sci. Tokyo U. (IA) {\bf 34} (1987), 223--243;
{\it IV}, J. Fac. Sci. Tokyo U. (IA) {\bf 36} (1989), 297--328.}
\item{[\xms]}
{P. Martin and H. Saleur, {\it On an Algebraic Approach to Higher Dimensional
Statistical Mechanics}, Comm. Math. Phys. {\bf 158} (1993), 155-190.}
\item{[\xmsa]}
{P. Martin and H. Saleur, {\it The blob algebra and the periodic
Temperley--Lieb algebra}, Lett. Math. Phys, {\bf 30} (1994) no. 3, 189--206.}
\item{[\xshi]}
{Shi Jian-yi, {\it The Kazhdan--Lusztig cells in the affine Weyl
group}, Lecture Notes in Mathematics, {\bf 1179}, Springer, New York, 1986.}
\item{[\xwes]}
{B.W. Westbury, {\it The representation theory of the Temperley--Lieb
Algebras}, \newline Math. Z. {\bf 219} (1995), 539--565.}

\vfill
\eject
\end